# ESTIMATION FOR ALMOST PERIODIC PROCESSES

By Keh-Shin Lii and Murray Rosenblatt

*University of California, Riverside and University of California, San Diego*

Processes with almost periodic covariance functions have spectral mass on lines parallel to the diagonal in the two-dimensional spectral plane. Methods have been given for estimation of spectral mass on the lines of spectral concentration if the locations of the lines are known. Here methods for estimating the intercepts of the lines of spectral concentration in the Gaussian case are given under appropriate conditions. The methods determine rates of convergence sufficiently fast as the sample size $n \to \infty$ so that the spectral estimation on the estimated lines can then proceed effectively. This task involves bounding the maximum of an interesting class of non-Gaussian possibly nonstationary processes.

**1. Introduction.** The main objective of this paper is to present a constructive method to determine the lines of support of spectra or, equivalently, the frequencies or periods of Gaussian harmonizable processes with almost periodic covariances under appropriate conditions. The processes considered provide an interesting set of random processes that are generally not transient and not stationary and for which Fourier methods of analysis are helpful and meaningful. The lines of support of spectra are lines parallel to the diagonal in the two-dimensional spectral plane. A number of papers have discussed spectral estimation in this context but they all assume knowledge of the lines of support of the spectra. The novelty of this paper is in the presentation of *constructive methods for estimating* the lines of support under appropriate conditions with a rate of convergence good enough to imply bias and covariance of spectral estimation using the estimated support lines with the same asymptotic properties as if the actual lines of support were known precisely.

Almost periodically correlated processes have been considered by Alekseev [1], Hurd [10], Gardner [5], Gladyshev [9], Gerr and Allen [7, 8], Dandawate and Giannakis [3], Leskow and Weron [14] and Lii and Rosenblatt [15]. Extensive application of periodically correlated (or cyclostationary) processes is noted in [6] and references therein. For the sake of simplicity and economy, we limit ourselves to a discrete time parameter process $\{X_n, n = \ldots, -1, 0, 1, \ldots\}$ with mean $EX_n \equiv 0$. The process $\{X_n\}$ with covariance function $r_{n,m} = EX_n \overline{X}_m$ is said to be *periodically correlated* with period $T$ if

$$r_{n+T, m+T} = r_{n,m}. \tag{1.1}$$

We shall call $\{X_n\}$ a process with *almost periodically correlated covariance function* if for each $(s, \tau)$, there are functions $\alpha_j(s, \tau)$ and real $\lambda_j(s, \tau)$ with $\sum_j |\alpha_j(s, \tau)| < \infty$ such that, for each $s, \tau$, $r_{s+t, \tau+t} = \sum_j \alpha_j(s, \tau) e^{i \lambda_j(s, \tau) t}$ with the convergence uniform in $s, \tau$. We shall call a sequence $\{a(n)\}$ almost periodic if there are sequences $\{\alpha_j\}$, real $\{\beta_j\}$ such that $a(n) = \sum_j \alpha_j e^{i \beta_j n}$ with $\sum_j |\alpha_j| < \infty$.

The process $\{X_n\}$ is harmonizable (in the sense of Loève) if it has a Fourier representation in mean square,

$$X_n = \int_{-\pi}^{\pi} e^{in\lambda} \, dZ(\lambda), \tag{1.2}$$

with the random spectral function having a covariance

$$\text{cov}(Z(\lambda), Z(\mu)) = F(\lambda, \mu),$$

$$r(n, m) = \int_{-\pi}^{\pi} \int_{-\pi}^{\pi} e^{in\lambda - im\mu} \, dF(\lambda, \mu), \qquad \int_{-\pi}^{\pi} \int_{-\pi}^{\pi} |dF(\lambda, \mu)| < \infty \tag{1.3}$$

[16, 20]. Almost periodicity of the covariance function of a harmonizable process implies that $F$ has mass (complex) on at most a countable number of diagonal lines

$$\lambda = \mu + b, \qquad b = b_j, \qquad j = \ldots, -1, 0, 1, \ldots.$$

If the process is real-valued (as we assume to be the case) and $\lambda = \mu + b$ is a line of spectral support, then so is $\lambda = \mu - b$.

If the spectra on the lines of support of $F$ are absolutely continuous with respect to Lebesgue measure with spectral densities $f_b(\mu)$, the assumption that $\{X_n\}$ is real-valued implies that

$$f_b(\mu) = \bar{f}_{-b}(\mu + b) = \bar{f}_{-b}(-\mu) = f_b(-\mu - b). \tag{1.4}$$

Also,

(1.5) $$|f_b(\mu)|^2 \leq f_0(\mu) f_0(\mu + b).$$

Let us call $\mathcal{L}$ the collection of $b$'s corresponding to lines with nontrivial spectral support. In [3], the spectra $f_b(\mu)$ are estimated under appropriate conditions assuming that the elements of $\mathcal{L}$ are known.

We consider estimation of $b \in \mathcal{L}$ under the assumption that $\mathcal{L}$ is finite and $\int_{-\pi}^{\pi} f_b(\lambda)\, d\lambda \neq 0$ and see whether the final estimates of $f_b(\mu)$ based on these estimated elements of $\mathcal{L}$ produce an effect.

At this point, a few additional remarks are made on related papers. Tian [18] considers time domain estimation of the integral period of a sequence with periodically correlated covariance function observed at integral time. Hurd and Gerr [11] and Gerr and Allen [7] use the coherence between spectral components to test for a specific period of a periodically correlated covariance function against the null hypothesis of stationarity. They used the integrated coherence along diagonal spectral support lines. Lund et al. [17] use an averaged squared coherence statistic to test existence of periodic correlation against stationarity. The case of a continuous parameter process $\{X(\cdot)\}$ is considered in [4] and [12].

Before the main discussion we describe a simple but interesting class of processes with almost periodic covariance function. Additional examples can be found in [6] and [20]. Set

(1.6) $$X_n = \sum_{j=1}^{m} \cos(w_j n) Y_{n-\tau_j}^{(j)}, \qquad 0 \leq w_j \leq \pi,$$

where the processes $\{Y_n^{(j)}\}$ are jointly stationary with means zero and the cross-spectral density of $Y^{(j)}(\cdot)$ and $Y^{(k)}(\cdot)$ is $f_{jk}(\cdot)$. This commonly used model in communication is a mixture of amplitude modulated (AM) signals with carrier frequency $w_j$ modulated by the signal $Y_n^{(j)}$ with possible relative delay $\tau_j$. If the stationary processes $\{Y_n^{(j)}\}$ have the spectral representation

(1.7) $$Y_n^{(j)} = \int_{-\pi}^{\pi} e^{in\lambda} \, dZ_j(\lambda),$$

the process $\{X_n\}$ is given by $X_n = \int_{-\pi}^{\pi} e^{in\lambda}\, dZ(\lambda)$ with

(1.8) $$dZ(\eta) = \tfrac{1}{2} \sum_j [e^{-i\tau_j(\eta - w_j)}\, dZ_j(\eta - w_j) + e^{-i\tau_j(\eta + w_j)}\, dZ_j(\eta + w_j)].$$

Since the vector of processes $\{Z_j(\lambda)\}$ has orthogonal increments, the spectral mass of the process $\{X_n\}$ is given by

$$(1.9)\ E\,dZ(\lambda)\overline{dZ(\mu)} = \begin{cases} \frac{1}{4}\sum_j \{f_{jj}(\lambda - w_j) + f_{jj}(\lambda + w_j)\}\,d\lambda, & \text{on } \lambda = \mu, \\ \frac{1}{4}\Bigg\{\sum_{j,k} e^{-i(\tau_j - \tau_k)(\lambda - w_j)} f_{jk}(\lambda - w_j) \\ \qquad + \sum_{j,k}{}' e^{-i(\tau_j - \tau_k)(\lambda + w_j)} f_{jk}(\lambda + w_j)\Bigg\}\,d\lambda, \\ \quad \text{where the first sum is over } (j,k) \text{ such that} \\ \quad \lambda - \mu = w_j \pm w_k = \eta \neq 0 \\ \quad \text{and the second primed sum over } (j,k) \text{ such that} \\ \quad \lambda - \mu = -w_j \pm w_k = \eta \neq 0, \\ 0, & \text{otherwise.} \end{cases}$$

The estimation method proposed in this paper can be used to estimate the frequencies $w_j$, the spectral densities of the $Y^{(j)}$'s and the $\tau_j$'s under certain assumptions. We note that the remarks for almost periodic sequences imply that a moving average process

$$X_n = \sum_{j=0}^{m} a_j(n)\xi_{n-j}, \qquad E\xi_n \equiv 0,$$

with the sequences (in $n$) $\{a_j(n)\}, j = 0, 1, \ldots, m$, almost periodic and the process $\{\xi_n\}$ a white noise process, is a process with almost periodic covariance function. That is also true of the autoregressive scheme $X_n - a_n X_{n-1} = \xi_n$ if $\{a_n = \sum_j \alpha_j e^{i\beta_j n}\}$ is an almost periodic sequence with $\sum_j |\alpha_j| < 1$ (see also [13]).

We note that in meteorological, astronomical and engineering applications questions of periodic or almost periodic covariances are often of interest (see [17] and [20] for examples). The organization of the rest of the paper is as follows. In the next section notation and assumptions are given. The main result on the estimation of intercepts of almost periodic spectral support lines using an integrated periodogram [see (2.6)] is presented after Corollary 4. The estimation method is based on Theorem 4. Theorems 5 and 7 show that the asymptotic bias and variance of spectral estimates on the estimated lines are the same as those when the support lines are known. Section 3 contains a simulated example to demonstrate the effectiveness of the methods. The Appendix contains most of the proofs of the results.

**2. Results.** In this section we introduce notation and assumptions, and state the main results. The principal thread of the derivation will be indicated. However, the detailed proofs of the results will be given in the Appendix. By

$$(2.1) \qquad W = Z \bmod 2\pi,$$

it is to be understood that $-\pi < W \leq \pi$ and if $\{u\}$ is the integer $\ell$ with $-1/2 < u - \ell \leq 1/2$, then $(Z \bmod 2\pi) = Z - \{Z/(2\pi)\}2\pi$. Consider the finite Fourier transform

$$(2.2) \qquad F_n(\lambda) = \sum_{t=-n/2}^{n/2} X_t e^{-it\lambda}$$

and the periodogram

$$(2.3) \qquad I_n(\lambda, \mu) = \frac{1}{2\pi(n+1)} F_n(\lambda) \overline{F_n(\mu)}.$$

Here we assume that $n$ is even and $t$ is summed over integers. Let $K(\eta)$ be a continuous nonnegative symmetric weight function with finite support and such that $\int_{-\infty}^{\infty} K(\eta) \, d\eta = 1$. Set $K_n(\eta) = b_n^{-1} K(b_n^{-1}\eta)$ with $b_n \downarrow 0$ and $nb_n \to \infty$ as $n \to \infty$. A plausible estimate of $f_w(\eta)$ is given by

$$(2.4) \qquad \hat{f}_w(\eta) = \int_{-\pi}^{\pi} I_n(\mu + w, \mu) K_n(\mu - \eta) \, d\mu.$$

ASSUMPTION 1. $\{X_k\}$, $EX_k = 0$, is a harmonizable Gaussian sequence with almost periodic covariance function. The set $\mathcal{L}$ of $b$'s with nontrivial support is finite, and the $f_b(\mu), b \in \mathcal{L}$, are continuously differentiable on the curve $\lambda = (b + \mu) \bmod 2\pi$ on the torus of points $(\lambda, \mu) \in [-\pi, \pi]^2$, with $(-\pi, \mu)$ identified with $(\pi, \mu)$ and $(\lambda, -\pi)$ with $(\lambda, \pi)$.

Now

$$(2.5) \qquad \begin{aligned} r_{n,m} &= EX_n X_m = \sum_{b \in \mathcal{L}} \int_{-\pi}^{\pi} e^{in(\mu+b)} e^{-im\mu} f_b(\mu) \, d\mu = \sum_{b \in \mathcal{L}} e^{inb} c_b(n-m), \\ c_b(n) &:= \int_{-\pi}^{\pi} e^{in\lambda} f_b(\lambda) \, d\lambda. \end{aligned}$$

ASSUMPTION 2. There is a function $c(j) \geq 0$ on the integers with $|c_b(j)| \leq c(j)$ for all $b \in \mathcal{L}$ and $\sum_j c(j) < \infty$. Unless otherwise noted, the summation over $j$ is on the integers.

ASSUMPTION 3. The $f_b(\mu), b \in \mathcal{L}$, are such that $\int_{-\pi}^{\pi} f_b(\mu) \, d\mu \neq 0$.

In general, it is possible that there may be spectral lines $b$ with $\int_{-\pi}^{\pi} f_b(\mu) \, d\mu = 0$ (see [11], page 347). In that case let $\mathcal{L}$ consist only of those $b$ with $\int_{-\pi}^{\pi} f_b(\mu) \, d\mu \neq 0$ and our arguments will be valid for the restricted $\mathcal{L}$. The local extrema of the absolute value of the integrated periodogram

$$(2.6) \qquad U_n(b) := \int_{-\pi}^{\pi} I_n(b + \mu, \mu) \, d\mu = \frac{1}{n+1} \sum_{k=-n/2}^{n/2} X_k^2 e^{-ibk}$$

*will be used to estimate the $b \in \mathcal{L}$ under Assumption* 3. *The fact that* (2.6) *is a periodogram computed for $X_k^2$ is due to the uniform weighting in integrating $I_n(b + \mu, \mu)$. With nonuniform weighting, one would get more complicated quadratic forms that might be useful if $\int_{-\pi}^{\pi} f_b(\mu)\, d\mu = 0$.*

Theorems 1 and 2 give estimates for $EI_n(\lambda, \mu)$ and $\operatorname{cov}(I_n(\mu+w, \mu), I_n(\mu'+w', \mu'))$ and it is *crucial* that they are given with error terms uniform in the $w$'s and $\mu$'s. Let

$$y(b, w) = (b - w) \bmod 2\pi. \tag{2.7}$$

It will also be convenient to make use of the following expressions for $b, b' \in \mathcal{L}$ when $\lambda = \mu + w, \lambda' = \mu' + w'$:

$$\begin{aligned}
y(1) &= (\lambda' + b - \lambda) \bmod 2\pi, \\
y(2) &= (-\mu' + \mu + b') \bmod 2\pi, \\
y(3) &= (-\mu' + b - \lambda) \bmod 2\pi, \\
y(4) &= (\lambda' + \mu + b') \bmod 2\pi.
\end{aligned} \tag{2.8}$$

Also let

$$\operatorname{sinc}(y) = \sin\left(\frac{n+1}{2} y\right) \Big/ \left(\frac{n+1}{2} y\right). \tag{2.9}$$

THEOREM 1. *If $\lambda = \mu + w$,*

$$EI_n(\lambda, \mu) = \sum_{b \in \mathcal{L}} f_b(\mu) \frac{\sin(((n+1)/2) y(b, w))}{((n+1)/2) y(b, w)} + O\left(\frac{\log n}{n}\right) \tag{2.10}$$

*with uniformity in $w, \mu$ when Assumption* 1 *is satisfied.*

THEOREM 2. *Assumption* 1 *and $\lambda = \mu + w, \lambda' = \mu' + w'$ imply that*

$$\begin{aligned}
\operatorname{cov}&(I_n(\mu+w, \mu), I_n(\mu'+w', \mu')) \\
&= \left\{ \sum_{b \in \mathcal{L}} f_b(\mu' + w') \operatorname{sinc}(y(1)) + O\left(\frac{\log n}{n}\right) \right\} \\
&\quad \times \left\{ \sum_{b' \in \mathcal{L}} f_{b'}(-\mu') \operatorname{sinc}(y(2)) + O\left(\frac{\log n}{n}\right) \right\} \\
&\quad + \left\{ \sum_{b \in \mathcal{L}} f_b(-\mu') \operatorname{sinc}(y(3)) + O\left(\frac{\log n}{n}\right) \right\} \\
&\quad \times \left\{ \sum_{b' \in \mathcal{L}} f_{b'}(\mu' + w') \operatorname{sinc}(y(4)) + O\left(\frac{\log n}{n}\right) \right\}.
\end{aligned} \tag{2.11}$$

ALMOST PERIODICITY AND ESTIMATION 7

Theorem 2 almost immediately implies Corollary 1.

COROLLARY 1. *Assumption 1 implies that*

$$\mathrm{cov}(\hat{f}_w(\eta), \hat{f}_{w'}(\eta'))$$

$$= \sum_{b,b'\in\mathcal{L}} \int_{-\pi}^{\pi}\int_{-\pi}^{\pi} f_b(\mu'+w')f_{b'}(-\mu')\,\mathrm{sinc}(y(1))\,\mathrm{sinc}(y(2))K_n(\mu-\eta)$$

$$\times K_n(\mu'-\eta')\,d\mu\,d\mu'$$

(2.12)

$$+ \sum_{b,b'\in\mathcal{L}} \int_{-\pi}^{\pi}\int_{-\pi}^{\pi} f_b(-\mu')f_{b'}(\mu'+w')\,\mathrm{sinc}(y(3))\,\mathrm{sinc}(y(4))$$

$$\times K_n(\mu-\eta)K_n(\mu'-\eta')\,d\mu\,d\mu'$$

$$+ O\left(\frac{\log n}{n}\right).$$

*Further,*

(2.13) $$\mathrm{cov}(\hat{f}_w(\eta), \hat{f}_{w'}(\eta')) = O\left(\frac{1}{nb_n}\right) + O\left(\frac{\log n}{n}\right).$$

Bounds on $\sup_w |U_n(w) - EU_n(w)|$ will be obtained by making use of Theorem 3. Let

(2.14) $$V_n(w) = (n+1)\{U_n(w) - EU_n(w)\}.$$

Cumulants of $V_n(w)$ will be estimated and an argument akin to Brillinger [2] will be used to gauge $\sup_w |V_n(w)|$.

THEOREM 3. *Let $\{X_k\}$ be a Gaussian sequence with almost periodic covariance function satisfying Assumptions 1 and 2. Then*

(2.15) $$\limsup_{n\to\infty} \sup_w |V_n(w)|/(n\log n)^{1/2} \leq 2^{5/2} q\left(\sum_j c(j)\right),$$

*with probability 1 where $q$ is the number of elements in $\mathcal{L}$.*

The arguments used to prove Theorem 3 will imply the following two corollaries.

COROLLARY 2. *Let $c > 0$. Then under the conditions of Theorem 3, there is a constant $C > 0$ such that*

(2.16) $$P\left[\sup_w |U_n(w) - EU_n(w)| > c\right] \leq \exp[-Cn^{1/2}(\log n)^{1/2}].$$

COROLLARY 3. *If $c > 0$ is fixed but large,*

$$(2.17) \quad P\left[\sup_w |U_n(w) - EU_n(w)| > 2cn^{-1/2}(\log n)^{1/2}\right] \leq \exp\left\{-\frac{c}{2}\log n\right\}.$$

An additional corollary indicates that $\sup_w |V_n(w)|$ diverges as $n \to \infty$.

COROLLARY 4. *Assume that $f_0(\cdot)$ is not identically zero. Then under Assumptions 1–3, there are a countable number of distinct $w_i$ in the neighborhood of zero such that, for each finite subcollection $\mathcal{C}$,*

$$(2.18) \qquad \frac{1}{\sqrt{n+1}} V_n(w_i), \qquad w_i \in \mathcal{C}$$

*are jointly asymptotically normal and independent with variances each greater than a constant $\zeta > 0$. Further,*

$$(2.19) \qquad 1 = o\left[\sup_w |V_n(w)|\right]$$

*in probability as $n \to \infty$.*

We now describe the estimation procedure. *Fix a given level $\delta > 0$. Consider the number of $b \in \mathcal{L}$ with*

$$(2.20) \qquad \left|\int_{-\pi}^{\pi} f_b(\mu)\, d\mu\right| > \delta.$$

*The following procedure will yield estimates $w_n(b)$ of $b \in \mathcal{L}$ satisfying (2.20) which are such that $n|w_n(b) - b| < a_n = n^{-1/5}$ except for a set of probability less than $O(\log n/n^2)$ if $n$ is sufficiently large (i.e., as $n \to \infty$). Divide the range (it can be taken as $-\pi \leq w \leq \pi$) of $w$ into blocks of length $c(n)/n$ [with $\lim_{n\to\infty} c(n) = \infty$ and $c(n) = o(n)$ as $n \to \infty$, e.g., $c(n) = \log n$ or $c(n) = n^{1/2}$]... and compute the global maximum of $|\int_{-\pi}^{\pi} I_n(\mu + w, \mu)\, d\mu|$ in each block. Throw out the blocks with maxima that fall less than or equal to $\delta + (\frac{\log n}{n})^{1/4}$. Also, if there are two maxima located within $c(n)/n$ of each other, throw out the block with the smaller maximum. As $n \to \infty$, the locations $w_n(b)$ at which the remaining maxima occur satisfying (2.20) are such that $n|w_n(b) - b| < a_n = n^{-1/5}$ except for a set of probability less than $O(\log n/n^2)$ as $n \to \infty$.*

We sketch the argument for this. Using (2.10), it is seen that

$$EU_n(w) = \sum_{b \in \mathcal{L}} \frac{\sin((n+1)/2(b-w))}{(n+1)/2(b-w)} \int_{-\pi}^{\pi} f_b(\mu)\, d\mu + O\left(\frac{\log n}{n}\right).$$

The maximum of

$$(2.21) \qquad \left|\frac{\sin((n+1)/2)(b-w)}{((n+1)/2)(b-w)}\right|\left|\int_{-\pi}^{\pi} f_b(\mu)\, d\mu\right|$$

is clearly at $w = b$ given Assumption 1. The contribution from the other terms in $|EU_n(w)|$ with $b' \neq b, b' \in \mathcal{L}$ is asymptotically $O(\frac{1}{n})$ in a neighborhood $\{w : |b - w| < \min_{b' \neq b, b' \in \mathcal{L}} \frac{1}{2}|b' - b|\}$ of $b \in \mathcal{L}$. First note that, uniformly in $w$,

$$(2.22) \qquad ||U_n(w)| - |EU_n(w)|| \leq C_1 \left(\frac{\log n}{n}\right)^{1/2}$$

except on a set of probability less than $O(n^{-2})$ by Corollary 3. We shall show that the local maximum $w_n(b)$ converging to $b$ must be such that $|n(w_n(b) - b)| < n^{-1/5}$ except for a set of probability less than $O(\log n/n^2)$. First note that

$$\left|\left|\int_{-\pi}^{\pi} I_n(\mu + w, \mu)\, d\mu\right| - \left|\int_{-\pi}^{\pi} f_b(\mu)\, d\mu \frac{\sin((n+1)/2)(b-w)}{((n+1)/2)(b-w)}\right|\right|$$
$$\leq C_2 \left(\frac{\log n}{n}\right)^{1/2}$$

in a fixed nontrivial interval $I(b, \varepsilon) = \{w : |b - w| < \varepsilon\}, 0 < \varepsilon < \frac{1}{2}$, containing $b$ but excluding the other points $b' \in \mathcal{L}$ except for a set of probability $O(\log n/n^2)$. The value of $|\int I_n(\mu + w, \mu)\, d\mu|$ at $w = b$ is at least

$$①:= \left|\int_{-\pi}^{\pi} f_b(\mu)\, d\mu\right| - C_2 \left(\frac{\log n}{n}\right)^{1/2}.$$

However, the value outside $|n(w - b)| < n^{-1/5}$, but inside $I(b, \varepsilon)$, of $|\int I_n(\mu + w, \mu)\, d\mu|$, is at most

$$②:= \left|\int_{-\pi}^{\pi} f_b(\mu)\, d\mu\right| \frac{\sin n^{-1/5}}{n^{-1/5}} + C_2 \left(\frac{\log n}{n}\right)^{1/2},$$

and since

$$\frac{\sin n^{-1/5}}{n^{-1/5}} = 1 - \frac{n^{-2/5}}{3!} + O(n^{-4/5}),$$

① is larger than ② for $n$ sufficiently large. Consequently, we have the following theorem.

THEOREM 4.  *Let $\delta > 0$ and Assumptions 1–3 hold. There are then estimates $w_n(b)$ of values $b \in \mathcal{L}$ with $|\int_{-\pi}^{\pi} f_b(\mu)\, d\mu| > \delta$ such that $n|w_n(b) - b| < n^{-1/5}$ except for a set of probability $O(\log n/n^2)$ if $n$ is sufficiently large.*

A similar but more detailed argument yields a better estimate for $|w_n(b) - b|$.

COROLLARY 5. *Under the assumptions of Theorem 4, $(n+1)(b-w_n(b)) = O(n^{-1/4}(\log n)^{1/4})$ and $[(n+1)(b-w_n(b))]^2 = O(\sup_w |U_n(w) - EU_n(w)|)$ except for a set of probability $O(\log n/n^2)$.*

The following results show that if one uses these estimates $w_n(b)$ of $b$ to estimate the spectral densities at $b$, the bias and variance are asymptotically as good as for the estimates with $b \in \mathcal{L}$ given.

THEOREM 5. *If $|\int_{-\pi}^{\pi} f_b(\mu)\,d\mu| > \delta$ and Assumptions 1–3 hold,*

$$E\hat{f}_{w_n(b)}(\eta)$$
$$= E\left[\frac{\sin((n+1)/2)y(b,w_n(b))}{((n+1)/2)y(b,w_n(b))} \int_{-\pi}^{\pi} f_b(\mu)K_n(\mu-\eta)\,d\mu\right] + O\left(\frac{\log n}{n}\right)$$

*and the right-hand side is*

$$f_b(\mu) + \frac{b_n^2}{2}f_b''(\mu)\int_{-\infty}^{\infty} x^2 K(x)\,dx + O(b_n^2) + O(n^{-1/2}(\log n)^{1/2}) + O\left(\frac{1}{nb_n}\right).$$

Corollaries 2 and 3 are used in the derivation of this result.

THEOREM 6. *Let $\delta > 0$ and Assumptions 1–3 hold. If*

$$\left|\int_{-\pi}^{\pi} f_b(\mu)\,d\mu\right| > \delta, \qquad \left|\int_{-\pi}^{\pi} f_{b'}(\mu)\,d\mu\right| > \delta$$

*for $b, b' \in \mathcal{L}$, then*

$$\operatorname{cov}(\hat{f}_{w_n(b)}(\eta), \hat{f}_{w_n(b')}(\eta')) = O\left(\frac{1}{nb_n}\right) + O\left(\frac{\log n}{n}\right)$$

*as $n \to \infty, b_n \downarrow 0, nb_n \to \infty$.*

THEOREM 7. *If $|\int_{-\pi}^{\pi} f_b(\mu)\,d\mu| > \delta$ and Assumptions 1–3 hold with $f_0(\cdot)$ not identically zero, setting $\gamma(x) = 1$ if $x = 0 \mod 2\pi$ and $0$ otherwise, then*

$$\operatorname{var}(\hat{f}_{w_n(b)}(\eta))$$
$$= \frac{2\pi}{nb_n}\left[f_0(\eta+b)f_0(-\eta)\right.$$
$$\left.+ \sum_{b' \in \mathcal{L}} \gamma(-2\eta+b'-b)|f_{b'}(-\eta)|^2\right]\int_{-\infty}^{\infty} K^2(x)\,dx + O\left(\frac{\log n}{n}\right).$$

A comparison of Theorems 5 and 7 with Proposition 5.1 and Corollary 5.1 in [15] indicate that the bias and variance of the spectral estimates using estimates $w_n(b)$ of $b \in \mathcal{L}$ have the same behavior asymptotically as spectral estimates assuming complete knowledge of $b \in \mathcal{L}$.

**3. Examples.** In this section we give a simulated example. Consider the amplitude modulated (AM) signal model in communication

$$(3.1) \quad y(t) = \cos(\omega_1 t)Z_1(t) + \cos(\omega_2 t)Z_2(t), \qquad EZ_i(t) = 0, \qquad i = 1, 2,$$

where $Z_1(t)$ and $Z_2(t)$ are independent stationary processes with spectral densities $f_1(\lambda)$ and $f_2(\lambda)$ at frequency $\lambda$, respectively. Then $\{y(t)\}$ is harmonizable with almost periodic covariance function and spectral support lines on $\lambda = \mu$, $\lambda = \mu \pm 2\omega_1$ and $\lambda = \mu \pm 2\omega_2$. We consider (3.1) with $\omega_1 = \frac{\pi}{4\sqrt{3}}, \omega_2 = \frac{\pi}{3\sqrt{2}}$, and $Z_1(t) = \varepsilon_t + 0.5\varepsilon_{t-1}, Z_2(t) = \eta_t + 0.3\eta_{t-1}$, where $\varepsilon_t$ and $\eta_t$ are independent standard Gaussian sequences. In this case, $y(t) = \int e^{it\lambda} dZ_y(\lambda)$ with $dZ_y(\lambda) = \frac{1}{2}[dz_1(\lambda \mp \omega_1) + dz_2(\lambda \mp \omega_2)]$ and $Z_i(t) = \int e^{it\lambda} dz_i(\lambda)$ for $i = 1, 2$. Therefore, the spectrum of $\{y(t)\}$ is given by

$$\begin{aligned}E[dZ_y(\lambda)\overline{dZ_y(\mu)}] \\ = \tfrac{1}{4}\{(f_1(\lambda \mp \omega_1) + f_2(\lambda \mp \omega_2))\delta_{\lambda,\mu} \\ + f_1(\lambda \mp \omega_1)\delta_{\lambda,\mu \pm 2\omega_1} + f_2(\lambda \mp \omega_2)\delta_{\lambda,\mu \pm 2\omega_2}\}.\end{aligned}$$

Given a realization of $\{y(t)\}$, we are interested in estimating $\omega_1$ and $\omega_2$ and the spectral densities on the lines given by $\lambda = \mu \pm 2\omega_1$ and $\lambda = \mu \pm 2\omega_2$.

As an example, a simulated realization of $\{y(t)\}$ for $t = -512$ to $512$ is generated. A stretch of $\{Z_1(t)\}$ and $\{y(t)\}$ are given in Figure 1. The periodogram of $\{y(t)\}$ is given in Figure 1 also. There is no apparent information on $\omega_1$ and $\omega_2$ in the periodogram. The $U_n(b)$ statistic of $\{y(t)\}$ is computed for a range of $b$'s. For $b = b_j = \frac{2\pi j}{1025}$ with $j = 1, 2, \ldots, 512$, $U_n(b)$ and $E(U_n(b))$ are given in Figure 2(a)–(f).

Two spikes are apparent at frequencies near $2\omega_1$ and $2\omega_2$. The largest two $b$'s in $|U_n(b)|$ occurred at $\hat{\omega}_1 = 0.45406$ and $\hat{\omega}_2 = 0.73938$. When the $U_n(b)$ are calculated at a more refined grid with $b = b_j = \frac{2\pi j}{4097}$, $j = 1, 2, \ldots, 2048$, the largest $|U_n(b)|$ occurred at $\hat{\omega}_1 = 0.45329$ and $\hat{\omega}_2 = 0.74015$. These latter values are closer to the true values $\omega_1 = 0.45345$ and $\omega_2 = 0.74048$. These are plotted in Figure 2(g) and (h). The scale in figures (g) and (h) differs from that in figures (e) and (f) since there is leakage at frequencies very close to zero because the computation at (e) and (f) is at a more refined set of frequencies. The zero frequency corresponds to the diagonal line which is the support line of the "power spectrum" of the process. Therefore, the value of $|U_n(0)|$ is relatively large. It is this mass that leaks.

A large value of $|U_n(b)|$ suggests that $\lambda = \mu + b$ is a spectral support line. Figure 3 gives the plots of spectral estimates on these estimated support lines. The dashed line gives the theoretical spectral density if the intercept $b$ is known exactly. In these cases the imaginary part is zero. For the diagonal line $\lambda = \mu$, there is no need to estimate the intercept. The estimated spectral densities on each line when the intercepts are exactly known are given by the solid lines. The dotted lines are spectral estimates when the intercepts are estimated with a crude search, while the dash-dot lines are spectral estimates when the intercepts are estimated by the more refined search.

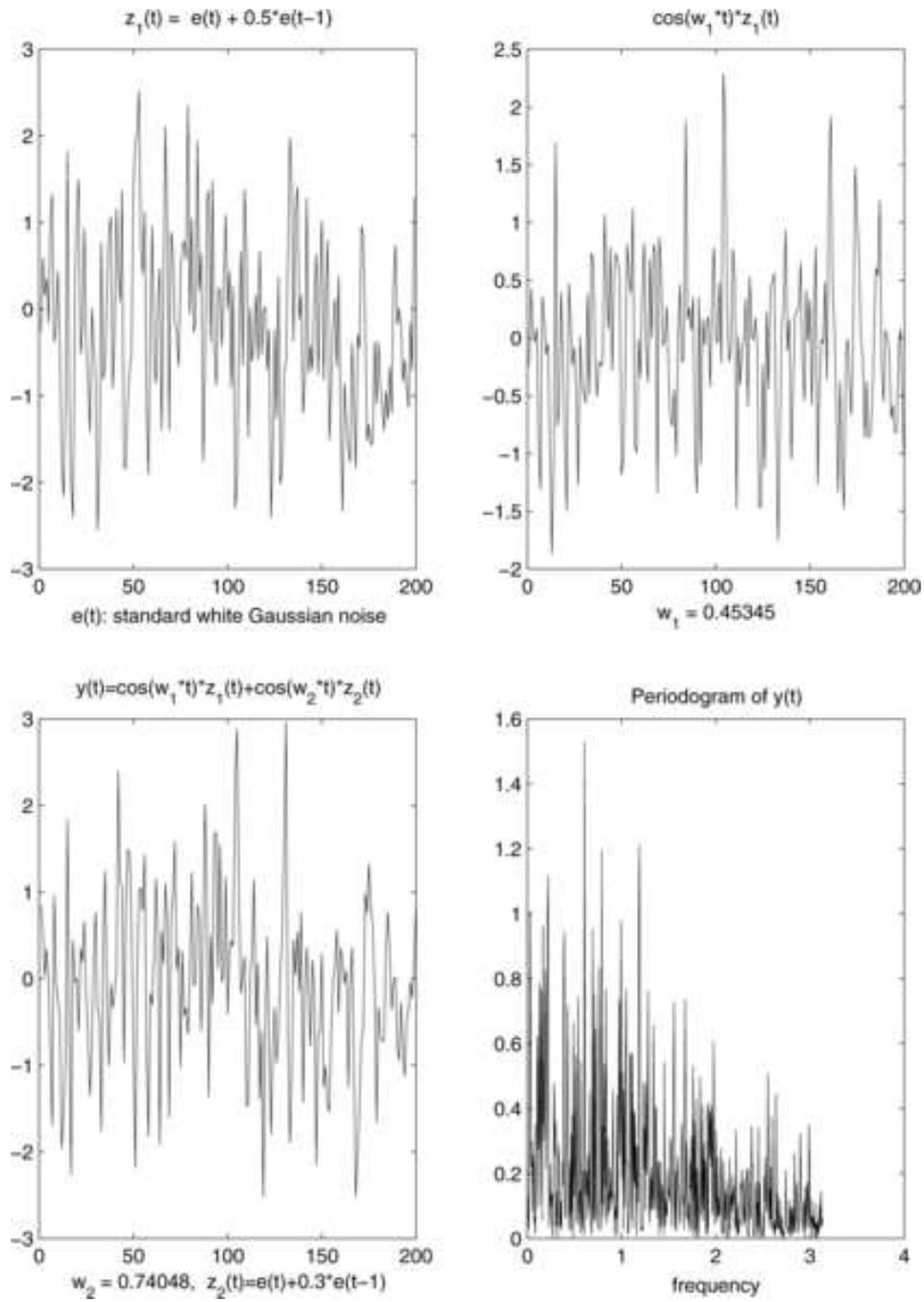

FIG. 1. *Simulated time series with $z_1(t), z_2(t)$ independent stationary time series.*

When the sample size of $\{y(t)\}$ is increased to $t = -2048$ to $2048$, we have more resolution as expected.

These graphs show that intercepts of these support lines or frequencies or "periods" of the covariance function can be effectively estimated using $U_n(b)$. These intercepts can be estimated sufficiently accurately to give accurate estimates of the spectral densities on the spectral support lines.

## APPENDIX

For the proofs, the following will be useful. Let

$$\text{(A.1)} \quad D_n(x) = \frac{\sin((n+1)/2)x}{\sin(x/2)}, \qquad \widetilde{D}_n(x) = \frac{\sin((n+1)/2)x}{x/2}.$$

LEMMA A.1. *If $|y| \leq \pi$,*

$$\text{(A.2)} \quad \int_{-\pi}^{\pi} D_n(x+y)D_n(-x)\,dx = \int_{-\pi}^{\pi} \widetilde{D}_n(x+y)\widetilde{D}_n(-x)\,dx + O(\log n),$$

*with $O$ uniform in $y$.*

The argument given in [15] for Lemma 7.2 holds for this case.

PROOF OF THEOREM 1. Given $X_t, t = -n/2, \ldots, n/2$, then

$$\text{(A.3)} \quad EI_n(\lambda, \mu) = \frac{1}{2\pi(n+1)} \sum_{b \in \mathcal{L}} \int_{-\pi}^{\pi} D_n(v+b-\mu-w)D_n(-v+\mu)f_b(v)\,dv$$

$$= \sum_{b \in \mathcal{L}} G_1,$$

with $G_1 = \frac{1}{2\pi(n+1)} \int_{-\pi}^{\pi} D_n(v - \mu + y(b,w)) D_n(-v+\mu) f_b(v)\,dv$. The continuous differentiability of $f_b$ implies that

$$\text{(A.4)} \quad f_b(v) = f_b(\mu) + (v-\mu)f_b'(v^*),$$

with $v^*$ between $\mu$ and $v$. Then $G_1 = G_2 + G_3$ with [denoting $y(b,w)$ by $y$]

$$G_2 = \frac{1}{2\pi(n+1)} \int_{-\pi}^{\pi} D_n(v-\mu+y)D_n(-v+\mu)(v-\mu)f_b'(v^*)\,dv,$$

$$G_3 = \frac{1}{2\pi(n+1)} f_b(\mu) \int_{-\pi}^{\pi} D_n(x+y)D_n(-x)\,dx.$$

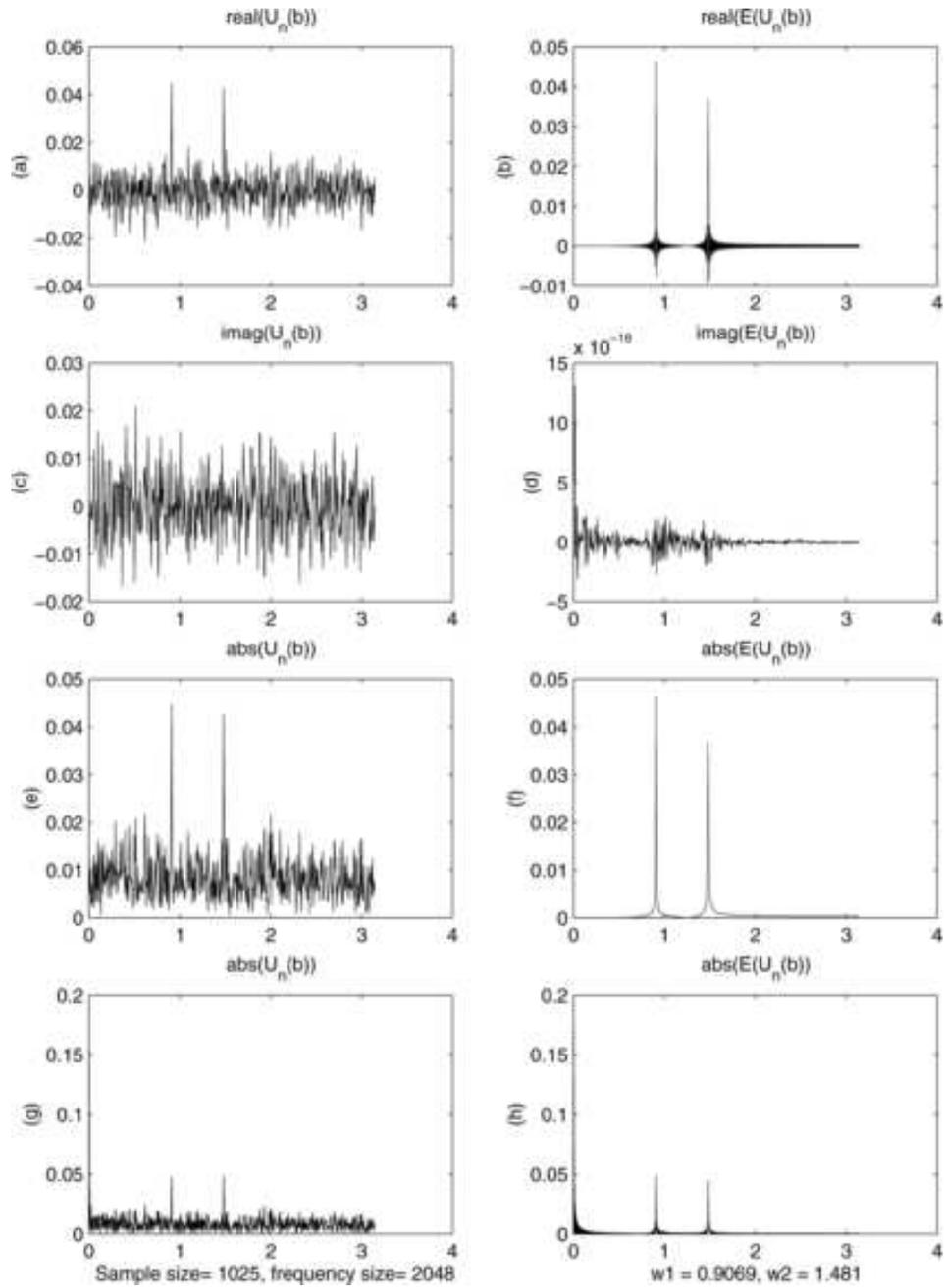

Fig. 2. $U_n(b)$ statistics of $\{y(t)\}$. Frequencies are computed for increments of $2\pi/1025$ for (a)–(f) and for increments of $2\pi/4097$ for (g)–(h).

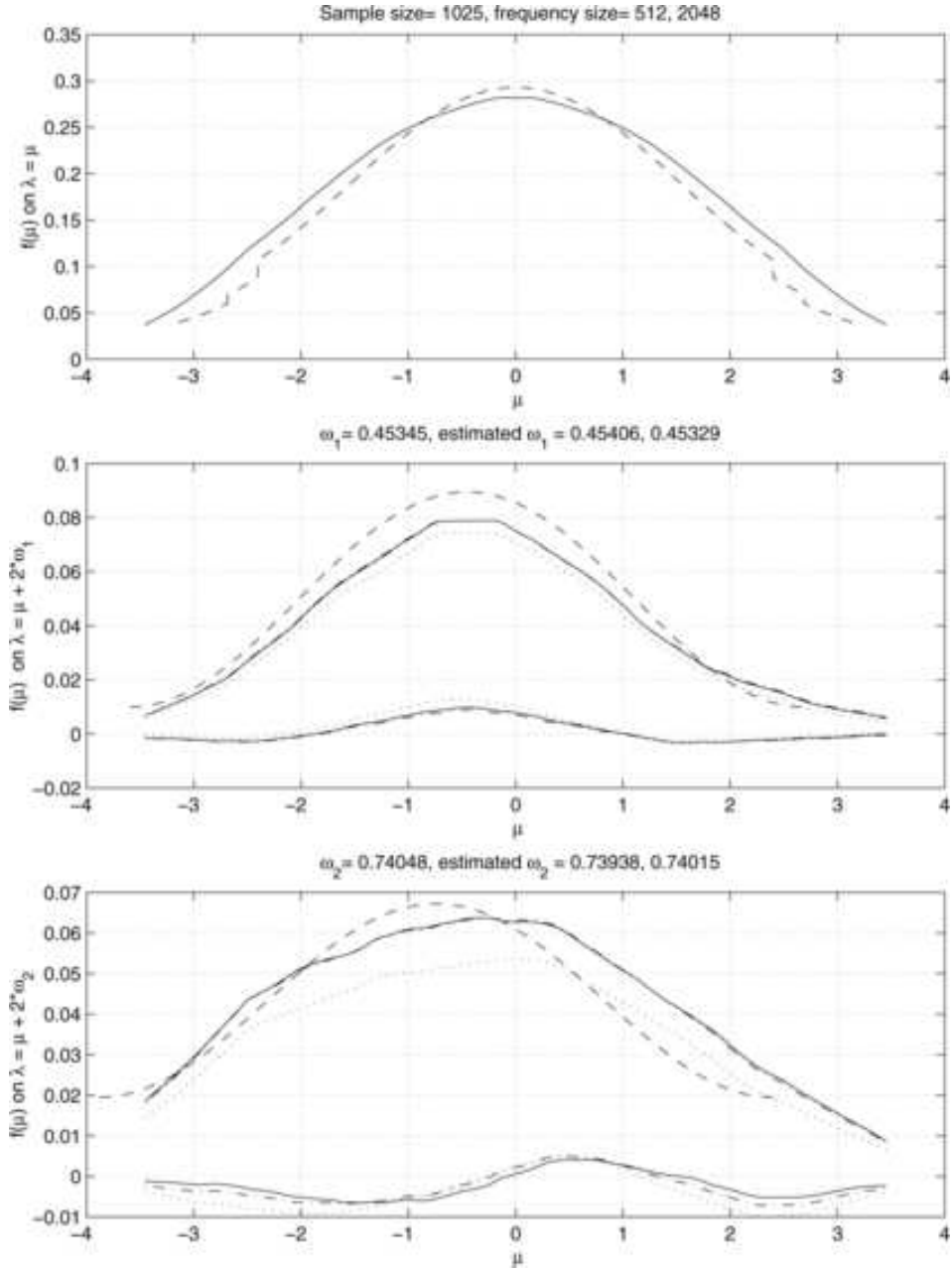

Fig. 3. *Spectral density estimates at different intercepts; the theoretical part is given by the dashed - - - line, estimated part with intercept known is given by the solid line; estimated part with intercept estimated at lower resolution is given by the dotted · · · line; estimated part with intercept estimated at higher resolution is given by dash-dot - · - line. Both real and imaginary parts are given. Theoretical imaginary parts are all zero and they are not plotted.*

It is clear that $|G_2| = O(\frac{\log n}{n})$. Also,

$$G_3 = f_b(\mu) \frac{1}{2\pi(n+1)} \int_{-\pi}^{\pi} \widetilde{D}_n(x+y)\widetilde{D}_n(-x)\,dx + O\left(\frac{\log n}{n}\right)$$

using Lemma A.1. To simplify the notation, we define

$$\left\{\sum_{i=1}^{I} \int_{A_i}\right\} f(x)\,dx = \sum_{i=1}^{I} \int_{A_i} f(x)\,dx.$$

Now

$$G_4 \equiv \frac{1}{2\pi(n+1)} \int_{-\pi}^{\pi} \widetilde{D}_n(x+y)\widetilde{D}_n(-x)\,dx$$

(A.5)
$$= \frac{1}{2\pi(n+1)}$$
$$\times \left\{\int_{-\infty}^{\infty} - \int_{|z|>(n+1)\pi}\right\} \frac{\sin(1/2)(z+(n+1)y)}{(z+(n+1)y)/2(n+1)} \frac{\sin(1/2)z}{z/2(n+1)} \frac{dz}{n+1}$$
$$= G_5 - G_6,$$

with $G_5$ and $G_6$ the two integrals of the previous line. Notice that

$$|G_6| \leq \frac{1}{2\pi}\left\{\int_{(n+1)\pi+1\geq|z|\geq(n+1)\pi} + \int_{2(n+1)\pi\geq|z|\geq(n+1)\pi+1} + \int_{|z|\geq 2(n+1)\pi}\right\}$$
$$\left|\frac{\sin(1/2)(z+(n+1)y)}{(z+(n+1)y)/2} \frac{\sin(1/2)z}{z/2}\right| dz,$$

with the integrand less than or equal to $1/((n+1)\pi)$ in the first integral, less than or equal to $\frac{1}{(n+1)\pi}(z+(n+1)y)^{-1}$ in the second integral and less than $2/z^2$ in the third integral (because $|y| \leq \pi$). These bounds are enough to insure that $G_6 = O(\frac{\log n}{n})$. In evaluating $G_5$, we note that

(A.6) $$\frac{\sin(1/2)(z+(n+1)y)}{(z+(n+1)y)/2} = \int_{-\infty}^{\infty} \mathbb{1}_{[-1/2,1/2]}(x)e^{ix(z+(n+1)y)}\,dx.$$

The Parseval relation implies that $G_5 = \text{sinc}(y)$. Applying the above bounds for the $G_i$'s proves the theorem. □

PROOF OF THEOREM 2. By a direct calculation, we have four expressions of the same type as that analyzed in Theorem 1. Using (2.8) and (2.10), formula (2.11) is obtained. □

PROOF OF COROLLARY 1. The first equation follows directly from Theorem 2. The order of magnitude statement is obtained by making use of the Schwarz inequality. □

PROOF OF THEOREM 3. First notice that, for distinct indices in the range $-n/2$ to $n/2$,

$$\text{(A.7)} \quad \text{cum}(X_{j_1}^2, X_{j_2}^2, \ldots, X_{j_s}^2) = 2^{s-1} \sum_{\nu} r_{j_1, j_{\alpha_2}} r_{j_{\alpha_2}, j_{\alpha_3}} \cdots r_{j_{\alpha_s}, j_1},$$

where $\nu$ is a listing of the $(s-1)!$ irreducible sequences $1 \to \alpha_2 \to \cdots \to \alpha_s \to 1$ with $\alpha_2, \ldots, \alpha_s$ distinct choices of $j_2, \ldots, j_s$. The case in which the indices are not distinct is obtained from (A.7) by making the appropriate identification of identical indices. By (2.5),

$$\text{(A.8)} \quad |r_{j_1, j_2}| \leq \sum_{b \in \mathcal{L}} |c_b(j_1 - j_2)|.$$

Since the $c_b$'s are uniformly bounded by a summable $c(\cdot) \geq 0$,

$$\text{(A.9)} \quad |r_{j_1, j_{\alpha_2}}, r_{j_{\alpha_2}, j_{\alpha_3}} \ldots r_{j_{\alpha_s}, j_1}| \leq q^s c(j_1 - j_{\alpha_2}) c(j_{\alpha_2} - j_{\alpha_3}) \cdots c(j_{\alpha_s} - j_1).$$

We shall carry out an argument for $\text{Re}\{V_n(b)\}$ alone since that for $Im\{V_n(w)\}$ is similar. Let $\sigma_n$ be the variance of $\text{Re}\{V_n(w)\}$. Then

$$\text{(A.10)} \quad \log E \exp[\alpha \text{Re}\{V_n(w)\}] - \alpha^2 \sigma_n / 2$$
$$\leq (n+1) \sum_{s=3}^{\infty} (2q)^s \alpha^s \left( \sum_{j=-\infty}^{\infty} c(j) \right)^s / s$$

since the $s$th cumulant of $V_n(b)$ is less than or equal to $(s-1)!(n+1)(2q)^s \times (\sum_j c(j))^s$. If $\alpha > 0$ is sufficiently small,

$$\text{(A.11)} \quad \log E \exp[\alpha \text{Re}\{V_n(b)\}] \leq (n+1) q^2 \left( \sum_j c(j) \right)^2 (1 + \epsilon) \alpha^2$$

since

$$0 \leq \sigma_n = \text{var}(\text{Re}\{V_n(w)\})$$
$$= \text{var}\left( \sum_{k=-n/2}^{n/2} X_k^2 \cos(kw) \right)$$
$$= 2 \sum_{j,k=-n/2}^{n/2} r_{j,k}^2 \cos(jw) \cos(kw)$$
$$< 2(n+1) q^2 \left( \sum_j c(j) \right)^2$$

if $\{X_k\}$ is not the trivial process that is identically zero. By Lemma 2.1 and Corollary 2.1 of [19] (see also [21], Chapter 10, Theorem 7.28),

$$\text{(A.12)} \quad \sup_w |\text{Re}\{V_n(w)\}| \leq \sup_j |\text{Re}\{V_n(u_j)\}| / (1 - 6\pi n R^{-1}),$$

with $u_j = 2\pi j/R, j = 0, 1, \ldots, R-1$ and $R \gg n$. The simple argument for this follows from the fact that $\max_\lambda |p'(\lambda)| \leq n \max_\lambda |p(\lambda)|$ for a trigonometric polynomial $p(\lambda) = \sum_{v=-n}^n \alpha_v \exp(iv\lambda)$ of order $n$ (Bernstein's inequality, page 11, [21], vol. 2). Then for $\lambda_j = 2\pi j/R, j = 0, 1, \ldots, R-1$,

$$\text{(A.13)} \qquad \max_{|\lambda| \leq \pi} |p(\lambda)| \leq \max_{|j| \leq R} |p(\lambda_j)| + \frac{2\pi}{R} \max_{|\lambda| \leq \pi} |p'(\lambda)|.$$

Equation (A.13) implies that

$$\left(1 - \frac{2n\pi}{R}\right) \max_{|\lambda| \leq \pi} |p(\lambda)| \leq \max_{|j| \leq R} |p(\lambda_j)|$$

and (A.12) follows. Then

$$\text{(A.14)} \qquad \begin{aligned} & E \exp\left\{\alpha \sup_w |\operatorname{Re}\{V_n(w)\}|\right\} \\ & \leq E \exp\left\{\frac{\alpha \sup_j |\operatorname{Re}\{V_n(u_j)\}|}{(1 - 6\pi n R^{-1})}\right\} \\ & \leq \sum_j E \exp\left\{\frac{\alpha |\operatorname{Re}\{V_n(u_j)\}|}{(1 - 6\pi n R^{-1})}\right\} \\ & \leq 2 \exp\left\{\log R + \frac{(1+\epsilon)(n+1)q^2(\sum c(j))^2 \alpha^2}{[(1 - 6\pi n R^{-1})^2]}\right\}. \end{aligned}$$

Our objective is to show by an appropriate choice of $a$ and $R$ that

$$\text{(A.15)} \qquad P\left[\sup_w |\operatorname{Re}\{V_n(w)\}| \geq 2a\right] < C n^{-1-\delta}$$

for some constant $C$. By (A.14) and Chebyshev's inequality, the probability on the left-hand side of (A.15) is less than or equal to

$$\text{(A.16)} \qquad \exp\{-2a\alpha\} 2 \exp\left\{\log R + (1+\epsilon)(n+1)q^2 \right. \\ \left. \times \left(\sum_j c(j)\right)^2 \alpha^2 / [(1 - 6\pi n R^{-1})^2]\right\}.$$

Given $\epsilon, \delta > 0$, set

$$\text{(A.17)} \qquad \begin{aligned} a^2 &= (1+\epsilon)(2+\delta)q^2(n+1)\log(n+1)\left(\sum_j c(j)\right)^2, \\ R &= n \log n \end{aligned}$$

and

(A.18) $$\alpha = a(1 - 6\pi n R^{-1})^2 \bigg/ \bigg[(1+\epsilon)(n+1)q^2\bigg(\sum_j c(j)\bigg)^2\bigg].$$

On inserting these expressions for $a, \alpha$ and $R$, (A.16) becomes

(A.19) $$2\exp[-(1+\delta)\log(n+1)] \leq Cn^{-1-\delta}.$$

The bound $Cn^{-1-\delta}$ has been obtained for sufficiently large $n$ with an appropriate choice for the constant $C$. Thus,

(A.20) $$P\bigg[\sup_w |\mathrm{Re}\{V_n(w)\}|/(n\log n)^{1/2} \geq 2\{(1+\epsilon)(2+\delta)\}^{1/2}q\bigg(\sum_j c(j)\bigg)\bigg]$$
$$\leq Cn^{-1-\delta}.$$

By applying the Borel–Cantelli lemma to (A.20), one can see that, for $\eta > 0$,

(A.21) $$\limsup_{n\to\infty}\sup_w |\mathrm{Re}\{V_n(w)\}|/(n\log n)^{1/2} \leq 2^{3/2}q(1+\eta)\bigg(\sum_j c(j)\bigg).$$

Since the same result holds for $\sup_w |ImV_n(w)|$, we obtain

(A.22) $$\limsup_{n\to\infty}\sup_w |V_n(w)|/(n\log n)^{1/2} \leq 2^{5/2}(1+\eta)q\bigg(\sum_j c(j)\bigg). \qquad \Box$$

PROOF OF COROLLARY 2. Notice that the probability of (A.15) can be rewritten as $P[\sup_w |\mathrm{Re}\{U_n(w) - EU_n(w)\}| \geq \frac{2a}{n+1}]$. Let $c > 0$ and $a = cn/4, R = n\log n$ and $\alpha = n^{-1/2}(\log n)^{1/2}$. Using (A.16), we see that then there is a $C > 0$ such that

$$P\bigg[\sup |\mathrm{Re}\{U_n(b) - EU_n(b)\}| > \frac{c}{2}\bigg] \leq \exp(-Cn^{1/2}(\log n)^{1/2}).$$

The same type of inequality holds for $Im\{U_n(b) - EU_n(b)\}$. Thus, the corollary holds. $\Box$

PROOF OF COROLLARY 3. The arguments are similar but now set $a = cn^{1/2}(\log n)^{1/2}, R = n\log n$ and $\alpha = n^{-1/2}(\log n)^{1/2}$ with $c > 0$ fixed but large. $\Box$

PROOF OF COROLLARY 4. Trivially,

$$\sup_w |U_n(w) - EU_n(w)|^2 \geq \frac{1}{2\pi}\int_{-\pi}^{\pi} |U_n(w) - EU_n(w)|^2\, dw$$
$$= \sum_{j=-n/2}^{n/2} (X_j^2 - EX_j^2)^2 \frac{1}{n+1}.$$

The covariance of $\frac{1}{n+1}(U_n(w) - EU_n(w))$ and $\frac{1}{n+1}(U_n(\eta) - EU_n(\eta))$ is

$$= \frac{1}{(n+1)^2} \sum_{j,k=-n/2}^{n/2} \operatorname{cov}(X_j^2, X_k^2) e^{-ijw+ik\eta},$$

with $\operatorname{cov}(X_j^2, X_k^2) = 2r_{j,k}^2$ since $\{X_t\}$ is a Gaussian process. But then

$$(1) := \frac{1}{(n+1)^2} \sum_{j,k=-n/2}^{n/2} \left\{ \sum_{b \in \mathcal{L}} e^{ibj} c_b(j-k) \right\}^2 e^{-ijw+ik\eta}$$

$$= \frac{2}{(n+1)^2}$$

$$\times \sum_{b,b' \in \mathcal{L}} \int_{-\pi}^{\pi} \frac{\sin((n+1)/2)(b+b'-w+\alpha)}{\sin(1/2)(b+b'-w+\alpha)} \frac{\sin((n+1)/2)(\eta-\alpha)}{\sin(1/2)(\eta-\alpha)}$$

$$\times f_b * f_{b'}(\alpha) \, d\alpha,$$

where "$*$" denotes convolution. If $\eta = w - b - b'$, we essentially get a Fejér kernel centered at this point and so

$$(1) = \frac{4\pi}{n+1} \sum_{b,b' \in \mathcal{L}} \delta(\eta - w + b + b') f_b * f_{b'}(\eta) + O\left(\frac{1}{n^2}\right)$$

if there are a finite number of $b$ in $\mathcal{L}$. Notice that, for $\eta = w$, we get

$$(1) = \frac{2\pi}{n+1} \sum_{b \in \mathcal{L}} f_b * f_{-b}(\eta) + O\left(\frac{1}{n^2}\right).$$

We first show joint asymptotic normality of $U_n(w) - EU_n(w)$ for any finite number of $w$'s. Also, there are an infinite number of $w_i$ in a neighborhood of zero such that $w_i - w_j + b + b' \neq 0$ for any $b, b' \in \mathcal{L}$. For joint asymptotic normality, we have to deal with

$$(A.23) \qquad \sum_{j=1}^{m} [U_n(w_j) - EU_n(w_j)] \alpha_j.$$

The argument for estimates of the cumulants of $V_n(b)$ leading to Theorem 3 still holds for the $s$th cumulant of (A.23) with the additional factor $(\sum |\alpha_j|)^s$. The asymptotic variance is

$$(A.24) \quad \frac{4\pi}{n+1} \sum_{j,j'} \alpha_j \alpha_{j'} \sum_{b,b' \in \mathcal{L}} \delta(w_j - w_{j'} + b + b') f_b * f_{b'}(w_j) + O\left(\frac{1}{n^2}\right).$$

We now consider conditions for positivity of the coefficient of $\frac{4\pi}{n+1}$.

We trivially get asymptotic normality for $\sqrt{n} \sum_{j=1}^{m} \alpha_j [U_n(w_j) - EU_n(w_j)]$ but with a possibly singular distribution because the $s$th cumulants with $s > 2$ all tend to zero. The asymptotic variance of $U_n(w_j) - EU_n(w_j)$ is $\frac{4\pi}{n+1} \sum_{b \in \mathcal{L}} f_b * f_{-b}(w_j) + O(\frac{1}{n^2})$. Now $f_0 * f_0(0) = \int_{-\pi}^{\pi} f_0(w)^2 \, dw > 0$ if $f_0$ is not identically zero. Further, for $b \in \mathcal{L}$,

$$f_b * f_{-b}(0) = \int_{-\pi}^{\pi} f_b(w) f_{-b}(-w) \, dw = \int_{-\pi}^{\pi} |f_b(w)|^2 \, dw \geq 0.$$

Thus, $\sum_{b \in \mathcal{L}} f_b * f_{-b}(0) > 0$ and this implies $\sum_{b \in \mathcal{L}} f_b * f_{-b}(\eta) > \zeta > 0$ for $\eta$ in some neighborhood of zero. But there are a countable number of $w_i \neq 0$ in this neighborhood such that $w_i - w_j + b - b' \neq 0$ for any $b, b' \in \mathcal{L}$. Therefore, the $\sqrt{n}[U_n(w_j) - EU_n(w_j)]$ are asymptotically jointly normal and independent with variance greater than $\delta > 0$. Therefore,

$$1 = o\left[\sup_{w_i} |\sqrt{n}[U_n(w_i) - EU_n(w_i)]|\right] = o\left[\sup_w |\sqrt{n}[U_n(w) - EU_n(w)]|\right]$$

in probability. We have the corollary. □

PROOF OF COROLLARY 5. The proof is similar to that of Theorem 4 and is omitted. □

PROOF OF THEOREM 5. If $\lambda = \mu + w$, Theorem 1 gives

$$EI_n(\lambda, \mu) = \sum_{\substack{b \in \mathcal{L} \\ |y(b,w)| \leq \pi}} f_b(\mu) \frac{\sin(((n+1)/2)y(b,w))}{((n+1)/2)y(b,w)} + O\left(\frac{\log n}{n}\right).$$

Since

$$\hat{f}_{w_n}(\eta) = \int_{-\pi}^{\pi} I_n(\mu + w_n, \mu) K_n(\mu - \eta) \, du,$$

it follows that

$$E\hat{f}_{w_n(b)}(\eta) = \sum_{b' \in \mathcal{L}} E\left[\frac{\sin(((n+1)/2)y(b', w_n(b)))}{((n+1)/2)y(b', w_n(b))}\right]$$
$$\times \int_{-\pi}^{\pi} f_{b'}(\mu) K_n(\mu - \eta) \, d\mu + O\left(\frac{\log n}{n}\right),$$
$$E\left[\frac{\sin(((n+1)/2)y(b, w_n(b)))}{((n+1)/2)y(b, w_n(b))}\right] = (1) + (2),$$

where

$$(1) := E\left[\frac{\sin(((n+1)/2)y(b, w_n(b)))}{((n+1)/2)y(b, w_n(b))}; |n(b - w_n(b))| \leq kn^{-1/4}(\log n)^{1/4}\right],$$
$$(2) := E\left[\frac{\sin(((n+1)/2)y(b, w_n(b)))}{((n+1)/2)y(b, w_n(b))}; |n(b - w_n(b))| > kn^{-1/4}(\log n)^{1/4}\right],$$

with $E[f(X); A] := \int_A f(x) \, dP(x)$ where $P(\cdot)$ is the distribution of $X$. Now

$$(1) = E[1 - [n(b - w_n(b))]^2/3!$$
$$+ o([n(b - w_n(b))]^2); |n(b - w_n(b))| \leq kn^{-1/4}(\log n)^{1/4}].$$

By Corollary 3, $1 \geq (1) \geq (1 - k^2 n^{-1/2}(\log n)^{1/2})(1 - \frac{1}{n})$, while $|(2)| \leq \frac{1}{n}$ if $k$ is large enough. Therefore,

$$\left| E\left[\frac{\sin((n+1)/2)y(b, w_n(b))}{((n+1)/2)y(b, w_n(b))}\right] - 1 \right| \leq k^2 n^{-1/2}(\log n)^{1/2}.$$

Now if $b' \neq b$ consider

$$E\left[\frac{\sin((n+1)/2)y(b', w_n(b))}{((n+1)/2)y(b', w_n(b))}\right] = (1') + (2'),$$

where

$$(1') := E\left[\frac{\sin((n+1)/2)y(b', w_n(b))}{((n+1)/2)y(b', w_n(b))}; |n(b - w_n(b))| > c\right],$$

$$(2') := E\left[\frac{\sin((n+1)/2)y(b', w_n(b))}{((n+1)/2)y(b', w_n(b))}; |n(b - w_n(b))| \leq c\right].$$

By Corollary 2 $|(1')| \leq \exp(-Cn^{1/2}(\log n)^{1/2}) + O(n^{-2}), |(2')| \leq C/n$. Notice that $|(\frac{n+1}{2})(b - w_n(b))|^2 \leq k \sup |U_n(w) - EU_n(w)|$ except for a set of probability $O(n^{-2})$, which implies

$$P[|(n+1)(b - w_n(b))| \leq cn^{-1/2}(\log n)^{1/2}] \geq 1 - \exp[-c' \log n]$$

for some $c'$ and so $P[|(n+1)(b - w_n(b))| > cn^{-1/2}(\log n)^{1/2}] \leq \exp[-c' \log n]$. These estimates imply that, if $f_b$ is twice continuously differentiable, $E\hat{f}_{w_n(b)}(\eta) = f_b(\eta) + \frac{b_n^2}{2} f_b''(\eta) \int x^2 K(x) \, dx + o(b_n^2) + O(n^{-1/2}(\log n)^{1/2})$. □

PROOF OF THEOREM 6. First notice that

$$\operatorname{cov}(\hat{f}_{w_n(b)}(\eta), \hat{f}_{w_n(b')}(\eta'))$$
$$= E[E[\hat{f}_{w_n(b)}(\eta)\hat{f}_{w_n(b')}(\eta')|w_n(b) = w_n, w_n(b') = w_n']]$$
$$- E\hat{f}_{w_n(b)}(\eta)E\hat{f}_{w_n(b')}(\eta')$$
(A.25) $\quad = E[\operatorname{cov}(\hat{f}_{w_n(b)}(\eta), \hat{f}_{w_n(b')}(\eta')|w_n(b) = w_n, w_n(b') = w_n')]$
$$+ E[E(\hat{f}_{w_n(b)}(\eta)|w_n(b) = w_n, w_n(b') = w_n')$$
$$\times E(\hat{f}_{w_n(b')}(\eta')|w_n(b) = w_n, w_n(b') = w_n')]$$
$$- E\hat{f}_{w_n(b)}(\eta)E\hat{f}_{w_n(b')}(\eta').$$

The estimate for the first term on the right-hand side of (A.25) is obtained by making use of Corollary 1. The remaining terms on the right-hand side of (A.25) can be seen to be

$$\sum_{b'',b'''\in\mathcal{L}} \text{cov}\left[\frac{\sin((n+1)/2)y(b'',w_n(b))}{((n+1)/2)y(b'',w_n(b))}, \frac{\sin((n+1)/2)y(b''',w_n(b'))}{((n+1)/2)y(b''',w_n(b'))}\right]$$

$$\times \int_{-\pi}^{\pi} f_b(\mu) K_n(\mu - \eta)\, d\mu \int_{-\pi}^{\pi} f_{b'}(\mu) K_n(\mu - \eta')\, d\mu + O\left(\frac{\log n}{n}\right).$$

Of course,

$$|\text{cov}(\text{sinc}\, y(b'', w_n(b)), \text{sinc}\, y(b''', w_n(b')))|$$
$$\leq \{\text{var}(\text{sinc}\, y(b'', w_n(b)))\, \text{var}(\text{sinc}\, y(b''', w_n(b')))\}^{1/2}.$$

Notice that $\text{var}(\text{sinc}\, y(b', w_n(b)))$ is

$$E[(\text{sinc}\, y(b'', w_n(b)))^2] - (E[\text{sinc}\, y(b'', w_n(b))])^2$$
$$= E[(1 + \{(\text{sinc}\, y(b'', w_n(b))) - 1\})^2] - [E(1 + \{\text{sinc}\, y(b'', w_n(b)) - 1\})]^2$$
$$= E[\{\text{sinc}\, y(b'', w_n(b)) - 1\}^2] - (E[\{\text{sinc}\, y(b'', w_n(b)) - 1\}])^2.$$

One has to consider the case in which $b = b''$ and that in which $b \neq b''$. Also, one should note that the expression above is less than or equal to $E[\{\text{sinc}\, y(b'', w_n(b)) - 1\}^2]$. If $b = b''$,

$$E[\{\text{sinc}\, y(b, w_n(b)) - 1\}^2; |n(b - w_n(b))| \leq kn^{-1/4}(\log n)^{1/4}]$$
$$= E\left[\left(\left(\frac{((n+1)/2)(b - w_n(b))}{3!}\right)^2 + o\left(\frac{((n+1)/2)(b - w_n(b))}{3!}\right)^2\right)^2;\right.$$
$$\left. |n(b - w_n(b))| < kn^{-1/4}(\log n)^{1/4}\right]$$
$$\leq C\log n/n,$$

while, by Corollaries 3 and 4,

$$E[\{\text{sinc}\, y(b, w_n(b)) - 1\}^2; |n(b - w_n(b))| > kn^{-1/4}(\log n)^{1/4}] \leq c/n$$

if $k$ is large enough. If $b \neq b''$,

$$E[\{\text{sinc}\, y(b'', w_n(b)) - 1\}^2; |n(b - w_n(b))| \leq c] \leq c/n$$

because $w_n(b)$ is close to $b$ and $b'' - w_n(b)$ is large, while

$$E[\{\text{sinc}\, y(b'', w_n(b)) - 1\}^2; |n(b - w_n(b))| > c]$$
$$\leq e^{-Cn^{1/2}(\log n)^{1/2}} + O(n^{-2})$$

follows by Corollary 2 as $n \to \infty$. The proof of Theorem 6 is complete. □

Proof of Theorem 7.

$$\text{var}(\hat{f}_{w_n(b)}(\eta)) = E(\text{var}(\hat{f}_{w_n(b)}(\eta)|w_n(b) = w_n))$$
$$+ E[E(\hat{f}_{w_n(b)}(\eta)|w_n(b) = w_n)^2] \quad \text{(A.26)}$$
$$- \{E\hat{f}_{w_n(b)}(\eta)\}^2.$$

We first consider the first term on the right-hand side:

$$\text{var}(\hat{f}_{w_n(b)}(\eta)|w_n(b) = w_n)$$
$$= \int_{-\pi}^{\pi} \int_{-\pi}^{\pi} \left[ \left\{ \sum_{b \in \mathcal{L}} f_b(\mu' + w_n) \operatorname{sinc}(y(1)) + O\left(\frac{\log n}{n}\right) \right\} \right.$$
$$\times \left\{ \sum_{b' \in \mathcal{L}} f_{b'}(-\mu') \operatorname{sinc}(y(2)) + O\left(\frac{\log n}{n}\right) \right\}$$
$$+ \left\{ \sum_{b \in \mathcal{L}} f_b(-\mu') \operatorname{sinc}(y(3)) + O\left(\frac{\log n}{n}\right) \right\}$$
$$\left. \times \left\{ \sum_{b' \in \mathcal{L}} f_{b'}(\mu' + w_n) \operatorname{sinc}(y(4)) + O\left(\frac{\log n}{n}\right) \right\} \right]$$
$$\times K_n(\mu - \eta) K_n(\mu' - \eta) \, d\mu \, d\mu'.$$

Making use of the symmetry of $K$ about zero on taking the expectation, we obtain

$$\text{var}(\hat{f}_{w_n(b)}(\eta)|w_n(b) = w_n)$$
$$= \frac{2\pi}{nb_n} \left[ f_0(\eta + b) f_0(-\eta) \right.$$
$$\left. + \sum_{b' \in \mathcal{L}} \gamma(-2\eta + b' - b)|f_{b'}(-\eta)|^2 \right] \int_{-\infty}^{\infty} K^2(x) \, dx + O\left(\frac{\log n}{n}\right).$$

The term

$$E[E(\hat{f}_{w_n(b)}(\eta)|w_n(b) = w_n)]^2 - \{E\hat{f}_{w_n(b)}(\eta)\}^2$$

can be shown to be $O(\frac{\log n}{n})$ just as in the discussion of the corresponding term in the proof of Theorem 6. $\square$

Department of Statistics  
University of California, Riverside  
Riverside, California 92521  
USA  
E-mail: ksl@stat.ucr.edu

Department of Mathematics  
University of California, San Diego  
La Jolla, California 92093  
USA  
E-mail: mrosenblatt@ucsd.edu